**A population level analysis of the gender gap in mathematics: Results on over 13 million**

**children using the INVALSI dataset**

**Authors:** D. Giofrè[1*], C. Cornoldi[2], A. Martini[3], & E. Toffalini[2**]

**Affiliations:**

[1]Disfor, University of Genoa, Italy

[2]Department of General Psychology, University of Padua, Italy

[3]Invalsi, National institute for the Italian educational system, Rome, Italy

*David Giofrè, Disfor, University of Genoa, Italy, Corso Andrea Podestà, 2, david.giofre@unige.it

[**]Enrico Toffalini, University of Padua, Italy, Via Venezia 8, enrico.toffalini@unipd.it

**Acknowledgement:** The authors would like to thank the INVALSI for granting access to the data. The views and opinions expressed in this article are those of the authors and do not necessarily reflect the official policy or position of the INVALSI. The authors would also like to thank Ingrid Boedker for commenting on an early version of the paper. The authors are also grateful to the Editor (Dr. Richard Haier), the Action Editor (Dr. Thomas Coyle), and the reviewers for their comments and suggestions.

**Manuscript accepted: 19 June 2020.**




**Abstract**

Whether males outperform females in mathematics is still debated. Such a gender gap varies across countries, but the determinants of the differences are unclear and could be produced by heterogeneity in the instructional systems or cultures and may vary across school grades. To clarify this issue, we took advantage of the INVALSI dataset, that offered over 13 million observations covering one single instructional system (i.e., the Italian system) in grades 2, 5, and 8, in the period 2010-2018. Results showed that males outperformed females in mathematics (and vice versa in reading), with gaps widening from the 2nd through to the 8th grade. The gender gap in mathematics was larger in the richer northern Italian regions (also characterized by greater gender equality) than in southern regions. This was not explained by average performance or fully accounted for by economic factors. No such north-south difference of the gap emerged in reading. Results are discussed with reference to the literature showing that the gender gap varies across world regions.

**A population level analysis of the gender gap in mathematics: Results on over 13 million**

**children using the INVALSI dataset**

Female underrepresentation in mathematics, and more broadly in science, technology, engineering and mathematics (STEM) fields, has been repeatedly observed and has profound effects on society (Ceci & Williams, 2011; Ceci, Williams, & Barnett, 2009; Cheryan, Ziegler, Montoya, & Jiang, 2017; Stoet & Geary, 2018). In contrast, males, who routinely perform better on internationally-recognized measures of mathematical literacy, tend to display poorer performance in reading (see OECD, 2016); with a lower presence than girls at the high end of the distribution of reading literacy (Reilly, Neumann, & Andrews, 2019). Recent evidence also has shown that there are more males than females with relative advantages in mathematics (i.e., their math achievement larger than their reading achievement), which results in a bias toward STEM fields (Stoet & Geary, 2015, 2018). Such findings indicate a potential predisposition towards STEM subjects on the part of males, potentially explaining part of the gender bias evident in the fields (Ceci & Williams, 2010).

Several theories have been advanced to explain male-female differences in mathematics and reading, implicating a range of factors from biology to socioeconomics (for extensive reviews and theoretical proposals see Geary, 2010; Halpern, 2012). Primary biological factors include brain size and morphology, evolutionary pressures, and prenatal hormones (Ceci et al., 2009; Eals & Silverman, 1994; Finegan, Niccols, & Sitarenios, 1992; van der Linden, Dunkel, & Madison, 2017). In contrast, others have argued that biologically-based cognitive sex differences indirectly contribute to academic sex differences. For instance, the male advantage in visuospatial abilities facilitates the use of spatial strategies for solving some types of mathematics problems and through this contributes to the male advantage (Geary, 1996).



Additional theories implicate inequality in the social and educational system, in particular with reference to lower GDP, lower economic opportunities, and earlier socialization processes (Breda & Napp, 2019; Ceci et al., 2009; Eagly & Wood, 1999; Else-Quest, Hyde, & Linn, 2010) in the observed gender performance gap. The gender gap in mathematics performance, in particular, has been proposed to reflect broader societal inequalities (Breda, Jouini, & Napp, 2018): more egalitarian societies should have smaller gender gaps in math, as one aspect of social inequalities may be represented by gender equality itself. In contrast, Stoet and Geary (2018), using recent PISA data (OECD, 2016), showed that the gender gap in math performance might follow the opposite pattern (i.e., the male superiority in mathematics was greater in the most gender equal countries), and they proposed that the life-quality pressure in less gender-equal countries could promote girls' greater engagement with STEM fields.

Several other explanations have been proposed to explain the gender gap in STEM fields. Halpern and co-authors (Halpern et al. 2007) performed a review on the extant literature and concluded that gender differences in mathematics and science did not directly evolve, but seem to be indirectly related to differences in interests, which are in turn related to specific brain and cognitive systems. This observation is also supported by a meta-analysis by Su, Rounds, Armstrong (2009), showing that men preferred working with objects, while women preferred working with people, with large effect sizes. This meta-analysis also indicated that men showed stronger realistic and investigative interests, again with a large effect size, while women showed stronger artistic, social, and conventional interests, with effect sizes ranging from small to moderate. Finally, meta-analytic estimates indicated that men performed better than women on specific measures of engineering, science, and mathematics interests, with effect sizes ranging from large for the former and small in latter two. Based on these premises, the authors concluded



that interests may have a critical role in gendered occupational choices and gender disparity in the STEM fields. These findings have been supported by several other authors, including Ceci and co-authors (Ceci, Williams, Barnett, 2009), who argued that biological evidence was inconclusive, while sociocultural factors were more promising. In particular, the authors proposed that women's preferences, both free and constrained, represented the most powerful explanatory factor of gender differences (e.g., Ceci et al., 2009; see also Ceci & Williams, 2010; 2011).

Even though data on the general distribution are usually considered, there is also a large literature focusing on the tails. In their meta-analysis Hyde and Mertz (2009) found that gender differences in mathematics were smallest in samples of the general population, grew larger with increasingly selective samples, and were largest for highly selected samples at the top tail of the distribution. This led the authors to conclude that there is a more evident male advantage at the upper tails of the distribution of mathematics ability, that is, a larger gender gap when the most difficult math problems are concerned (Hyde & Mertz, 2009). This finding was replicated several times using different measures and samples (Halpern et al., 2007). In particular, Strand and co-authors found more males at the top 5% and 10% in quantitative skills, examining a very large national sample of children aged 11-12 years. Furthermore, males were also overrepresented in the lower tails in reading, both 5% and 10% (Strand, Deary, Smith, 2006). Data from the Study of Mathematically Precocious Youth (SMPY), an extremely large study on children with very high performance in mathematics, found that males outperformed females in mathematics, with small effect size, while virtually no differences were found on verbal abilities (Benbow & Stanley, 1980, 1983), with differences remaining relatively stable 20 years later (Benbow, Lubinski, Shea, & Eftekhari-Sanjani, 2000). This finding has also been found using



PISA data, showing that at the highest levels of math performance, the gap between males and females increases and the reading gap decreases, while for low performing students, the math gap disappears and the reading gap widens (Stoet & Geary, 2013).

In sum, there is a large body of evidence concerning gender differences in mathematics and reading, but researchers have described and explained them in different ways. The current study aims to add to the current literature by examining the same patterns found in previous studies but using regional differences across a single nation tested at several time points. This unique approach provides several advantages. First, educational systems, as well as societal pressures, teaching styles, and even topics covered in the teaching of mathematics, might largely differ across countries, and using a sample from the same population could overcome this problem. Furthermore, the actual number of recent international evaluations is rather limited (PISA, for example, occurs every three years, meaning that only three evaluations are available in a 9-year time span). This is not a problem *per se*, as it would be unlikely to see large changes in less than three years, but changes from one evaluation to the other, although based on very large samples, may occur due to specific contingencies (e.g., nature of the specific items used) and statistical artifacts (e.g., regression toward the mean), rather than reflecting real differences in the population of interest (see Hunt, 2011 on this issue). A meta-analytic approach, with random-effect models to account for the variability of the effects across surveys, might overcome this problem. For all these reasons, detailed national datasets with data from many consecutive years can give additional insight.

It should also be noted that results concerning a single education system, at the same time characterized by important regional differences, may offer other important information. Italy has large differences between the richer northern regions and the poorer southern regions, with



regions in central Italy representing an intermediate case, and the gap in social opportunities

between males and females is especially evident in the southern part of the country (Checchi &

Peragine, 2010). Felice (2011), describing the state of the literature on the north vs. south Italy

divide (also referred to in Italian as the *questione meridionale* [southern affair]), comments that

this divide has a very long history, from the end of the nineteenth century to the present day, and

is persistent over the long run. The causes of this divide are currently debated, with some

scholars holding the view that at the time of Unification (1861) northern regions were more

prosperous because of a better geographical position, more favorable natural endowments, and

higher human and social capital. At the moment, differences between northern and southern

Italian regions are quite striking (e.g., GDP in the northern Italian regions is about two times

higher compared to southern Italian regions, while the opposite is true for unemployment levels

which are extremely high in southern regions, compared to norther Italian ones; see Felice, 2011;

2012 for an extensive discussion).

A more recent analysis of the phenomenon argues that although the south as a whole may

have ranked somewhat lower than the center-north, it was misleading to consider southern Italy

as a uniform area (see Felice, 2011) and subtle differences should be considered at the regional

level. In fact, the way that Italian geographical areas are organized into regions offers an

appropriate unit, as Italian regions have separate governments, reflect different historical and

cultural backgrounds, were historically part of different states up until about 160 years ago, and

currently have different dialects (e.g., Neapolitan, which was the official language of the

Kingdom of the Two Sicilies), cuisines, artistic tradition and other differences. Analyses carried

out in several studies stressed the importance of examining the achievement gap at the regional



level rather than at the nationwide level or macro-areas (e.g., Cornoldi, Belacchi, Giofrè, Martini, & Tressoldi, 2010; Cornoldi, Giofrè, & Martini, 2013; Lynn, 2010).

Italy represents a unique opportunity for studying the gender gaps in achievement also because an extremely large dataset considering both mathematics and reading tasks is available for scrutiny. This dataset includes assessments of both mathematical and reading academic performance across different grades on the entire population over the last ten years. Surprisingly, this very large database has not yet been systematically analyzed to address this important issue. Gender differences in mathematics can therefore be studied with the advantage of having one of the largest samples ever studied in this field, including different years of assessment, different geographic areas and different grades, but within the same school system using the same curriculum (i.e., the same school curriculum is taught in every school in Italy) and the same age established for entering into primary school (the calendar year in which a child will turn 6). It is worth noting that the Italian Constitution, promulgated in 1948, established the presence of 20 regions. These regions are now relatively independent in some respects, such economic and financial, but not all areas (e.g., defense or education). In particular, the Italian school system is regulated at the national level. Schools in Italy depend on the Minister of Education in Rome, with programs and curricula established at a national level and taught in every single region. This national regulation also includes the request that all Italian students are every year administered the same Mathematical and Reading tests, with the consequence that a potential huge mass of information is collected.

For the present analyses we took advantage of this dataset considering an extremely large number of data points (over 13 million assessments in mathematics and reading) while relevant



social and economic moderating factors were retrieved from The Italian National Institute of Statistics, ISTAT (http://dati.istat.it).

Based on previous evidence on this topic we expected to find gender differences in mathematics, favoring males, and in reading, favoring females. We also expected to find differences between Italian regions with northern regions outperforming southern regions in both mathematics and reading. In addition, we also wanted to examine whether the two variables interacted, examining, for the specific case of mathematics, the so-called gender equality paradox, showing that the areas with higher levels of gender equality have also larger gaps in STEM fields (Stoet & Geary, 2018). Additional analyses were also performed to understand the presence of possible moderators of the gender gap, and to examine whether males are overrepresented in the upper and in the lower tails of the distribution.

## Methods

### Materials

Data from the INVALSI (Istituto Nazionale per la VALutazione del Sistema di Istruzione [National Institute for the Assessment of the Instruction System]) surveys on the Italian population ($2^{nd}$, $5^{th}$, and $8^{th}$ grades) for mathematics and reading from 2010 to 2018 were analyzed.

Each year, the INVALSI assesses the whole Italian child population in $2^{nd}$, $5^{th}$, and $8^{th}$ grades. In the Italian school system, children typically attend the $2^{nd}$ grade at 7 years of age; the $5^{th}$ grade at 10 years of age; and the $8^{th}$ grade at 13 years of age. The surveys are population-based, meaning that all Italian children in the same grade are tested. Therefore, the same cohorts of children are tested across different years at different grades. To ensure anonymity, however,



the INVALSI datasets are organized in such a way that it is not possible to associate the performances of any single child longitudinally.

The INVALSI tasks are inspired by their international counterparts such as PISA, TIMSS, and PIRLS. Specifically, the tasks consist of a series of around 30-45 items, tailored for each academic year. Evaluations occur in the child's class and typically last for two hours. Most items have a multiple-choice response format and the remaining have an open-choice response format (with only one possible correct answer). Tasks differ for each grade and year of survey. Each response is coded as either correct (1) or incorrect (0). For both mathematics and reading literacy, the items are created based on IRTs. Different areas are tested in mathematics, including algebra, generally presented in the form of realistic problems (e.g. requiring students to solve the following problem: "Iodine 131 halves its mass every 8 days due to radioactive decay. How long will it take for it to reduce from 2 grams to 0.250 grams?"); geometry (e.g. requiring students to calculate the area of a geometric figure given the length of its sides); and statistics and calculation (e.g. requiring students to understand a bar plot, calculate the number of possible combinations of a series of items). For reading literacy, the items are distinguished into two main areas: reading comprehension (involving both fictional and informational texts followed by multiple choices questions), and grammar (e.g. requiring students to tick the instances in which a given verb expresses a future tense). More information about the test, with additional examples, is available at the INVALSI website https://invalsi-areaprove.cineca.it/

In total, we obtained 13,347,861 observations for mathematics, and 13,279,322 observations for reading literacy. In addition to the total population, the INVALSI survey offers, for each year, data on a selected representative sub-sample. Data collection on this sample is supervised by members of the INVALSI committee themselves, therefore guaranteeing no



cheating or other biases. We further analyzed these data to assess whether the same conclusions were reached as for the total population. Observations for the representative samples (which differed in each year of survey) were collected on 658,271 children for mathematics, and 655,945 for reading literacy, for the period between 2010 and 2017 (for 2018, no such data were available to us).

**Data analysis**

The meta-analytic estimates for the gender gaps in mathematics and reading were computed on standardized differences in performance between males and females (positive values indicate males' advantage; negative values indicate females' advantage). We adopted a meta-analytic approach even though we had virtually the entire child population over the mentioned time course, because data can be still seen as a sample of a much wider set of possible observations. In fact, most variability (and with it the non-negligible uncertainty on the final estimates) must be attributed to differences in results across different years (as shown below), which is most likely due to differences between alternative versions of the tasks. Following the same approach, other estimates, based on meta-analytic estimates, were calculated, including correlation between mathematics and reading, and odds ratios between males and females on the tails of the distributions. The meta-analytic strategy for all analyses followed the guidelines suggested by Borenstein and colleagues (Borenstein, Hedges, Higgins, & Rothstein, 2009).

For each combination of grade and year of survey (corresponding to an average of about 490,000 children per combination) we used Rasch models to determine the individual performance parameters (Rasch scores) of each child. These correlated very strongly (about $r =$



.99) with those calculated by the INVALSI team themselves, and with the scores calculated as the simple proportions of correct responses for each child.

The effect of interest for the present study was the standardized difference (Cohen's *d*) between the performance of males and the performance of females. First, it was calculated separately for each Italian region in each year of survey and each grade, and its variance was calculated using the formula for the variance of standardized difference suggested by Borenstein et al. (2009). To simplify the interpretation, a meta-analytic model was fitted separately for each grade. Random-effect models were used in order to account for the variability of the effects across surveys conducted in different years (i.e., heterogeneity). Regions of Italy were entered as the fixed moderating factor of interest in the models, ordered from the northernmost to the southernmost regions along the main northwest-to-southeast axis that characterizes the Italian territory (Italy is roughly tilted by 45° rather than following a straight north-to-south direction; see the specific order in the vertical axis of any figure in the Supplemental material). According to the INVALSI system, 8 Italian regions are classified as "North", 4 regions as "Center", 8 regions as "South".

## Results

A preliminary analysis, based on meta-analytic estimates using the Fisher's z transformation (Borenstein et al., 2009), showed that mathematics and reading scores were strongly and similarly correlated across all grades: $2^{nd}$ grade, $r = .59$, 95% CI [.58, .59]; $5^{th}$ grade, r = .62, [.61, .63]; and $8^{th}$ grade, $r = .63$ [.63, .64]. These correlations varied negligibly across different regions and grades (all $\Delta r < .08$).



Results on gender gaps showed that male superiority in mathematics increased with grade and the same (with an opposite superiority of females) happened for reading. However, while the gender gap in reading performance was relatively homogeneous and increased steadily with grade across Italian regions, the gender gap in mathematics performance followed a sharp geographical gradient, with the most northern Italian regions presenting a larger gender gap than the most southern regions (Figure 1; see also Figure S1 and S2 in Supplemental material for details, and Figure S3, which shows that the results were virtually identical for the supervised and more controlled sub-sample of the INVALSI data).

Figure 1 about here

## Analysis of moderators

In-depth analyses were conducted on data for the mathematics performance of 8[th] graders, because this is the highest grade, thus the gender gaps are the largest. Unsurprisingly, heterogeneity explained nearly all variability of the effect, $I^2 = 98\%$, as the very large pool of observations meant that the sampling error was nearly null, and the remaining variability was due to the random effect of year of survey or moderators. Heterogeneity of the effects across years is clearly depicted in all figures in the supplemental material. The analysis of moderators, presented in the supplemental material, has an exploratory purpose.

***Region and geographical gradient***. An obvious moderator is the "region" factor, which alone explained 70% of the heterogeneity. However, an excellent proxy of it was the simple geographical gradient, observed when regions were assigned an integer number representing their sequential order from the northernmost to the southernmost one (see the ordered list in y axis of Figure S1, S2, Supplemental material). Such a gradient explained 57% of the



heterogeneity of the effect, only slightly less than in the case when regions were simply considered as a generic factor. The Bayesian Information Criterion (BIC), strongly favored this more parsimonious model, $\Delta BIC$ = -38.99.

*Baseline performance*. Average mathematical achievement was higher in the northern than in the southern regions, $d$ = .06 ($d$ = .15 in the supervised sub-sample). However, this was not a significant moderator of the gender gap, $\chi^2(1)$ = 2.10, $p$ = .15. See Figure S4 in supplemental material for a simultaneous representation of both average performance and gender gap by region. The average performance of females in mathematics was rather similar in northern (following the North-South direction the interested regions are Trentino-Alto-Adige through Emilia-Romagna according to the INVALSI classification system) vs. southern (Abruzzo to Sicily) regions, Cohen's $d$ = .01 95% CI [-.12, .13]. Conversely, some differences were observed in males, with higher performances in northern vs. southern regions, Cohen's $d$ = .15 [05, .25].

*Time*. Time (i.e., year treated as a quantitative variable) was a significant, albeit modest, moderator of the male-female gap in mathematics, $\chi^2(1)$ = 17.64, $p$ < .001, $\Delta R^2$ = 10% (on residual heterogeneity, i.e., over and above the effect of the geographical gradient), $B$ = -.007 [-.011, -.004], suggesting that over a period of 10 years there was an estimated reduction of $\Delta d$ = -.07 [-.11, -.04]. Unsurprisingly, the moderating effect of time interacted with the geographical gradient itself, $\chi^2(1)$ = 17.85, $p$ < .001, $\Delta R^2$ = 11%, which implies that the gender gap in mathematics is shrinking in the regions where it exists (i.e., northern and central regions in particular), but not in the southern regions, were it was already approaching zero.

*Economic and social factors*. Relevant and well-known differences among northern and southern Italian regions concerning socio-economic factors were analyzed. We therefore



collected data on a series of socioeconomic factors: i) related to the geographical gradient; and,

ii) potential moderators of the gender gap in mathematics. These data were retrieved and

elaborated from the Italian National Institute of Statistics, ISTAT: http://dati.istat.it.

*Per capita GDP*, measured year-by-year throughout the period of the INVALSI survey,

was strongly correlated with the geographical gradient (on average, Spearman's $\rho = .89$; see

Figure 2). However, per capita GDP explained only 44% of the heterogeneity in the gender gap

in mathematics (as an alternative to "region" as a factor or gradient), and the effect of the

geographical gradient, included in a second step, was explaining a portion of the variance over

and above GDP, $\chi^2(1) = 42.42$, $p < .001$, $\Delta R^2 = 23\%$.

Figure 2 about here

*Average employment/unemployment* data (2004-2018) were also considered (see Figure

S5 for further details). Both employment and unemployment ratios followed the geographical

gradient, Spearman's $\rho = -.92$, and $\rho = .91$, respectively, showing a clearly disadvantaged

condition of the southern regions. For both variables, females were at a disadvantage vis-à-vis

males, and this was clearer in southern regions. The Odds Ratio (OR) was used as a measure of

the relative disadvantage of females over males in employment and unemployment. In fact, only

the OR (female/male) for employment followed the geographical gradient, Spearman's $\rho = .80$,

whereas the OR (female/male) for unemployment did not, Spearman's $\rho = .04$. Therefore, the

disadvantage of females over males in employment was more evident in the southern than in the

northern regions, but the same was not true for unemployment. Such a finding was probably due

to the fact that many unemployed females in southern Italian regions are not registered as

unemployed.



Notably, the total employment and unemployment ratios, or the same ratios separately considered for males and females, were extremely correlated, with Spearman's $\rho$s between .96 and .99. Therefore, the total employment ratio was chosen as the candidate for moderator of the gender gap in mathematics, as it also represents unemployment ratio (or the same ratios for females). This moderator alone explained 65% of heterogeneity of the gender gap in mathematics (lower employment ratios, especially of females, are related with smaller disadvantage of females in mathematical achievement). This is nearly all of the heterogeneity previously explained by region and/or the geographical gradient; the latter was no longer relevant when added as another moderator after employment ratio, $\chi^2(1) = 1.53$, $p = .21$, $\Delta R^2 = 1\%$. Region as a factor was not significant either, $\chi^2(18) = 25.25$, $p = .12$, $\Delta R^2 = 15\%$.

*Participation in local politics by males and females* data were retrieved for the period between 2008 and 2017 (see Figure S6 for further details). The proportion of females in city, province, and regional councils moderately followed the geographical gradient, with relatively more female participation in northern regions, as compared to southern regions, Spearman's $\rho = .55$. This factor was also strongly linked with the population size and historical relevance of the regions (e.g., Lazio, Campania and Sicily in Central-Southern Italy have a relatively high female share), meaning that this might be related to the gender gap but to a lesser extent. In fact, when included as a moderator of the gender gap in mathematics, the female share in participation in politics reached statistical significance, but it only explained a modest portion of the heterogeneity, $\chi^2(1) = 39.47$, $p < .001$, $\Delta R^2 = 20\%$. Such a finding indicates that lower female involvement in politics is linked, if only to a modest extent, with reduced gender gap in mathematics. In fact, the geographical gradient, included in a subsequent step, still explained a



very large portion of the residual variance after accounting for this factor, $\chi^2(1) = 102.62$, $p <$ .001, $\Delta R^2 = 47\%$.

***Gender ratios at the tails of the distributions.*** We also analyzed the ratios of males vs. females at the lower and upper tails of the distributions. For simplicity of reporting, we focused only on the 8[th] grade, where the gender gaps were larger, and grouped Italian regions into three clusters: North, Center, and South.

As the INVALSI surveys are an instrument for measuring achievement at the population level, they allow us to reliably identify children up to the top 1%, but not the top 0.1% or more, as children with the highest performance reach a ceiling level and children with the lowest performance reach a floor level. Therefore, we focused on the top (and bottom) 1% and 5% of the population's achievement in mathematics and reading.

To estimate how males and females were represented at the tails of the distributions, we calculated odds ratios (ORs). These are the ratios between the odds of a male being at the higher (or lower) tail of the distribution, and the same odd for a female. These are roughly equivalent to the simple gender ratios often reported in the literature, but odds ratios are proper effect sizes, can be meta-analyzed, are independent from any potential unbalance between males and females in the sample, and should therefore be preferred (Cumming, 2012). For meta-analyzing odds ratios, we followed the suggestions by Borenstein et al. (2009).

Importantly, if achievement scores were normally distributed and homoscedastic (i.e., if they had the same variance in males and females), gender ratios at the tails could be determined *a priori* based on the standardized difference (i.e., Cohen's *d*) between males and females. To assess how observed data corresponded to data derived from these *a priori* expectations based on



normality, we simulated 10 million data that reproduced the gender gap observed in the

INVALSI data, and that were set to be normally distributed and homoscedastic.

Results (Table 1) show that males were overrepresented in the upper, but not in the lower

tail in mathematics. In both the upper and lower tails of mathematics performance, however,

males were slightly more represented than one would expect if scores were normally distributed

and homoscedastic (i.e., the observed OR and their CIs are higher than the "simulated" OR in the

lower tail), at least in the North and in the Center of the country. In reading this pattern does not

appear (see the discussion for an in-depth analysis of these results).

Table 1 about here

**Mean differences across reading and math**

Finally, we calculated a series of scores expressing relative advantage in mathematics

vis-à-vis overall achievement (i.e., mathematics minus the average between mathematics and

reading literacy). However, the standardized differences between males and females in this

adjusted score followed virtually the same path as the mathematics score alone. This was

arguably due to the mathematics scores following a clear geographical gradient, whereas reading

literacy scores did not (Figure S7).

**Discussion**

The issues concerning gender differences in mathematics and regional differences seem

of particular relevance today, also because success in mathematics could affect girls' interest in

scientific disciplines and, indirectly, the socioeconomic gender gap. This is particularly

important because mathematics and reading skills are associated with educational achievement,



higher grades, greater attained qualifications and higher adult socioeconomic status (e.g.,

McGee, Prior, Williams, Smart, & Sanson,2002; Ritchie & Bates, 2013).

Several previous studies on the gender gap in mathematics have used the same

international dataset (e.g., PISA) which is useful for comparing school systems (and gender

gaps) across countries. However, this approach does not allow us to disentangle how different

education systems and other co-occurring social differences across countries could affect the

gender gaps in mathematics. For this reason, the INVALSI dataset could offer new pieces of

evidence within one country, with the advantage of having the same education system in all

regions, but also large regional differences in other societal aspects. This dataset provides data on

the entire population, thus offering the opportunity of evaluating these aspects with an extremely

large dataset on several consecutive assessments. Using this dataset, we could therefore address

some of the previous problems encountered by using other sources of data (including PISA, for

example), and collect some very insightful data on the gender gap in mathematics.

Our results suggest that Italian males outperform females in mathematics, and the gender

gap in both mathematics and reading widens as a function of grade and level of performance.

Substantial regional differences also emerged in mathematics, with a larger gender gap in the

most northern regions, which has tended to diminish, albeit slightly, in recent years. This

regional difference was only found in mathematics, while the gender gap was consistent in all

Italian regions in reading (where, on the contrary, females outperform males, in every single

region and in every single assessment). The male-female gap did not seem to be related to the

average academic performance in mathematics and was only partially related with the *per capita*

GDP. We also calculated a series of scores expressing relative advantage in mathematics vis-à-

vis overall achievement (i.e., mathematics minus the average between mathematics and reading



literacy). The results showed that, also due to the fact that mathematics scores followed a clear geographical gradient, whereas reading literacy did not, the standardized differences between males and females in this adjusted score followed virtually the same path as the mathematics score alone.

Alternative explanations could be offered for these results. Results are also consistent with a bio-ecological explanation, meaning that when the environment is generally more supportive (e.g., in socio-economically more developed regions), the largest gender differences emerge (Bronfenbrenner & Ceci, 1994). However, socioeconomic factors seem crucial in producing the mathematical gender gap. In particular, life-quality pressures and larger gender inequality in southern regions (Checchi & Peragine, 2010) could promote females' engagement with mathematics, as suggested by international comparisons using PISA data (Stoet & Geary, 2018). As mentioned in the introduction, there is also literature indicating that interest is one of the strongest predictors of gender differences in mathematics and reading (Ceci et al., 2009; Ceci & Williams, 2010, 2011; Halpern et al., 2007; Su et al., 2009). It could be therefore argued that there's very little evolutionary selection for the types of mathematics assessed on these tests (e.g., PISA or INVALSI), and that the most developed regions allow for greater expression of individual interests and development of basic competencies, such as spatial abilities, that could influence performance in mathematics.

The present analyses suggest that the mathematical gender gap could vary across different sociocultural contexts (regions in this case), even within the same country. Factors contributing to such differences may be difficult to disentangle, also because they are mostly confounded. We suggest, however, that general socioeconomic factors could be crucially related with these differences. Level of unemployment seems to explain a large part of the variance in



the gender gap, in the sense that where the unemployment (including the female unemployment, but not only) is high, the gender gap in mathematics disappears. Unemployment is also related to gender inequality (an aspect for which reliable regional data could not be calculated), having been argued to be crucial in determining gender differences (Stoet & Geary, 2018), but seems to offer a more general explanation of the gender gap results. A possible interpretation of this result is that, when unemployment is high, males and females react differently to the difficult requests of mathematical instruction: the former with a lower engagement, while the latter with a higher engagement than in other regional areas with lower unemployment rates.

It also can be argued that in southern areas there is an underground economy, including tourism and low technology activities, that is attracting adolescent boys and young men and that does not require a strong education. Under these conditions, boys might disengage more from school than girls with the consequence of having lower scores across academic areas. This particular finding is in line with the results in mathematics, in which boys in southern Italian regions seem to be particularly impaired. However, the performance in reading seems not to confirm this observation, as in this case boys are impaired in reading to a similar extent in all Italian geographic regions. As for participation in local politics, this variable seems not to be a crucial determinant for explaining the gender gap. One possible explanation is that there is very little cross-region variance on this variable, with female participation levels lower than 25% in all Italian regions. Some insightful information could also be obtained from the analyses of extreme cases in the two tails. As for mathematics, there were more males than females among students with very high performances, but not among those with very low performances. To some degree, this may seem obvious as males had on average higher performances than females in mathematics. Interestingly, however, the proportion of males in both the upper and lower tail



was larger than could be expected under the assumptions of normality and homoscedasticity, suggesting that the male distribution is more variable. This pattern holds true for the central and northern regions, in which the gender gap in mathematics was more evident. This larger-than-expected (assuming normality and homoscedasticity) proportion of males in both tails is in line with results from other countries, such as the UK and the USA (e.g., Benbow, & Stanley, 1980; Strand et al., 2006). It is worth noting, however, that in these countries a larger proportion of males were found in the lower tail of the distribution, while this pattern was not evident in our data. This pattern of results shows that in Italy is the gender gap in mathematics is more pronounced in the north and center areas, while weaker in the south, not only with respect to the whole population but also, and to a greater extent, for the top 1% performance, with northern regions having an OR of 2.18 (i.e., there are more than two males for each female in the top performing students) and only 1.27 in southern Italian regions. One possible interpretation of the sex difference in mathematics across regions is that girls are more engaged in school and mathematics in the southern regions than are boys, but this effect should be investigated further.

As for reading, there was a large proportion of males in the bottom of this distribution, while the number of males was particularly low in top performing students. This was unsurprising, as males performed worse than females in reading. Also, this was rather homogeneous across geographical areas, which is consistent with the homogeneity of the gender gap in reading across all regions. The proportion of males in the lower tail of the reading distribution was slightly smaller than one could expect assuming normality, meaning that males are not overrepresented in the lower tail of performance, at least in reading.

To conclude, we have described mathematics and reading performance in an extremely large and national representative sample of children. The results for mathematics indicate that



gender differences tend to be higher in northern and central regions. This particular finding is also associated with an overrepresentation in the same regions of males in the two tails of the mathematical performance distribution, that is, of males performing extremely well and extremely poorly. Results for reading do not follow a geographical gradient and show that females outperform males in a similar way in every single region and at every level of performance, confirming that this phenomenon is rather general and uniformly distributed in the population. With this paper we have therefore described an important pattern of data and provided some speculation on the possible causes; however, we recognize that the causes of the gender gap in both mathematics and reading are very elusive and difficult to examine and cannot be fully addressed with the extant available data. Therefore, further research is needed to understand the cognitive and non-cognitive underpinnings of this gap.



**Data availability**

The data that support the findings of this study are available from INVALSI. Restrictions

apply to the availability of these data, which were used under license for this study. Data are

available https://invalsi-serviziostatistico.cineca.it/ with the permission of INVALSI.



# References


Benbow, C. P., & Stanley, J. C. (1980). Sex differences in mathematical ability: Fact or artifact?
    *Science, 210*(4475), 1262–1264. https://doi.org/10.1126/science.7434028

Benbow, C., & Stanley, J. (1983). Sex differences in mathematical reasoning ability: more facts.
    *Science, 222*(4627), 1029–1031. https://doi.org/10.1126/science.6648516

Benbow, C. P., Lubinski, D., Shea, D. L., & Eftekhari-Sanjani, H. (2000). Sex differences in
    mathematical reasoning ability at age 13: Their status 20 years later. *Psychological Science,
    11*(6), 474–480. https://doi.org/10.1111/1467-9280.00291

Borenstein, M., Hedges, L., Higgins, J., & Rothstein, H. (2009). *Introduction to meta-analysis*.
    Chichester, UK: John Wiley & Sons, Ltd.

Breda, T., Jouini, E., & Napp, C. (2018). Societal inequalities amplify gender gaps in math.
    *Science.Sciencemag.Org*, *359*, 1219–1220. Retrieved from
    https://science.sciencemag.org/content/359/6381/1219.short

Breda, T., & Napp, C. (2019). Girls' comparative advantage in reading can largely explain the
    gender gap in math-related fields. *Proceedings of the National Academy of Sciences*,
    201905779. https://doi.org/10.1073/pnas.1905779116

Bronfenbrenner, U., & Ceci, S. J. (1994). Nature-nurture reconceptualized in developmental
    perspective: a bioecological model. *Psychological Review*, *101*(4), 568–586.
    https://doi.org/10.1037/0033-295X.101.4.568

Ceci, S. J., & Williams, W. M. (2010). Sex differences in math-intensive fields. *Current
    Directions in Psychological Science*, *19*(5), 275–279.




https://doi.org/10.1177/0963721410383241

Ceci, S. J., & Williams, W. M. (2011). Understanding current causes of women's underrepresentation in science. *Proceedings of the National Academy of Sciences*, *108*(8), 3157–3162. https://doi.org/10.1073/pnas.1014871108

Ceci, S. J., Williams, W. M., & Barnett, S. M. (2009). Women's underrepresentation in science: Sociocultural and biological considerations. *Psychological Bulletin*, *135*(2), 218–261. https://doi.org/10.1037/a0014412

Checchi, D., & Peragine, V. (2010). Inequality of opportunity in Italy. *The Journal of Economic Inequality*, *8*, 429–450. https://doi.org/10.1007/s10888-009-9118-3

Cheryan, S., Ziegler, S. A., Montoya, A. K., & Jiang, L. (2017). Why are some STEM fields more gender balanced than others? *Psychological Bulletin*, *143*(1), 1–35. https://doi.org/10.1037/bul0000052

Cornoldi, C., Belacchi, C., Giofrè, D., Martini, A., & Tressoldi, P. (2010). The mean Southern Italian children IQ is not particularly low: A reply to R. Lynn (2010). *Intelligence*, *38*(5), 462–470. https://doi.org/10.1016/j.intell.2010.06.003

Cornoldi, C., Giofrè, D., & Martini, A. (2013). Problems in deriving Italian regional differences in intelligence from 2009 PISA data. *Intelligence*, *41*(1), 25–33. https://doi.org/10.1016/j.intell.2012.10.004

Cumming, G. (2012). *Understanding the new statistics: Effect sizes, confidence intervals, and meta-analysis*. Routledge.

Eagly, A. H., & Wood, W. (1999). The origins of sex differences in human behavior: Evolved dispositions versus social roles. *American Psychologist*, *54*(6), 408–423.



https://doi.org/10.1037/0003-066X.54.6.408

Eals, M., & Silverman, I. (1994). The Hunter-Gatherer theory of spatial sex differences: Proximate factors mediating the female advantage in recall of object arrays. *Ethology and Sociobiology*, *15*(2), 95–105. https://doi.org/10.1016/0162-3095(94)90020-5

Else-Quest, N. M., Hyde, J. S., & Linn, M. C. (2010). Cross-national patterns of gender differences in mathematics: A meta-analysis. *Psychological Bulletin*, *136*(1), 103–127. https://doi.org/10.1037/a0018053

Felice, E. (2011). Regional value added in Italy, 1891-2001, and the foundation of a long-term picture. *The Economic History Review, 64*(3), 929–950. https://doi.org/10.1111/j.1468-0289.2010.00568.x

Felice, E. (2012). Regional convergence in Italy, 1891–2001: testing human and social capital. *Cliometrica, 6*(3), 267–306. https://doi.org/10.1007/s11698-011-0076-1

Finegan, J.-A. K., Niccols, G. A., & Sitarenios, G. (1992). Relations between prenatal testosterone levels and cognitive abilities at 4 years. *Developmental Psychology*, *28*(6), 1075–1089. https://doi.org/10.1037/0012-1649.28.6.1075

Geary, D. C. (1996). Sexual selection and sex differences in mathematical abilities. *Behavioral and Brain Sciences, 19*(02), 229. https://doi.org/10.1017/S0140525X00042400

Geary, D. C. (2010). *Male, female: The evolution of human sex differences*, 2nd edition. Washington, DC: American Psychological Association.

Halpern, D. F. (2012). *Sex differences in cognitive abilities, 4th edition*. New York, NY: Psychology Press

Halpern, D. F., Benbow, C. P., Geary, D. C., Gur, R. C., Hyde, J. S., & Gernsbacher, M. A.



(2007). The science of sex differences in science and mathematics. *Psychological Science in the Public Interest, 8*(1), 1–51. https://doi.org/10.1111/j.1529-1006.2007.00032.x

Hunt, E. (2011). *Human Intelligence*. Cambridge University Press.

Hyde, J. S., & Mertz, J. E. (2009). Gender, culture, and mathematics performance. Proceedings of the National Academy of Sciences, 106(22), 8801–8807. https://doi.org/10.1073/pnas.0901265106

Lynn, R. (2010). In Italy, north–south differences in IQ predict differences in income, education, infant mortality, stature, and literacy. *Intelligence*, *38*(1), 93–100. https://doi.org/10.1016/J.INTELL.2009.07.004

McGee, R., Prior, M., Williams, S., Smart, D., & Sanson, A. (2002). The long-term significance of teacher-rated hyperactivity and reading ability in childhood: Findings from two longitudinal studies. *Journal of Child Psychology and Psychiatry, 43*(8), 1004–1017. https://doi.org/10.1111/1469-7610.00228

OECD. (2016). *PISA 2015 Results (Volume I)*. https://doi.org/10.1787/9789264266490-en

Reilly, D., Neumann, D. L., & Andrews, G. (2019). Gender differences in reading and writing achievement: Evidence from the National Assessment of Educational Progress (NAEP). *American Psychologist, 74*(4), 445–458. https://doi.org/10.1037/amp0000356

Stoet, G., & Geary, D. C. (2013). Sex differences in mathematics and reading achievement are inversely related: Within- and across-nation assessment of 10 years of PISA data. *PLoS ONE*, *8*(3), e57988. https://doi.org/10.1371/journal.pone.0057988

Stoet, G., & Geary, D. C. (2015). Sex differences in academic achievement are not related to political, economic, or social equality. *Intelligence, 48*, 137–151.



https://doi.org/10.1016/j.intell.2014.11.006

Stoet, G., & Geary, D. C. (2018). The gender-equality paradox in science, technology,

engineering, and mathematics education. *Psychological Science*, *29*(4), 581–593.

https://doi.org/10.1177/0956797617741719

Strand, S., Deary, I. J., & Smith, P. (2006). Sex differences in Cognitive Abilities Test scores: A

UK national picture. *British Journal of Educational Psychology, 76*(3), 463–480.

https://doi.org/10.1348/000709905X50906

Su, R., Rounds, J., & Armstrong, P. I. (2009). Men and things, women and people: A meta-

analysis of sex differences in interests. *Psychological Bulletin, 135*(6), 859–884.

https://doi.org/10.1037/a0017364

Ritchie, S. J., & Bates, T. C. (2013). Enduring links from childhood mathematics and reading

achievement to adult socioeconomic status. *Psychological Science, 24*(7), 1301–1308.

https://doi.org/10.1177/0956797612466268

van der Linden, D., Dunkel, C. S., & Madison, G. (2017). Sex differences in brain size and

general intelligence (g). *Intelligence*, *63*, 78–88. https://doi.org/10.1016/j.intell.2017.04.007



**A.** Males advantage over females in mathematics performance

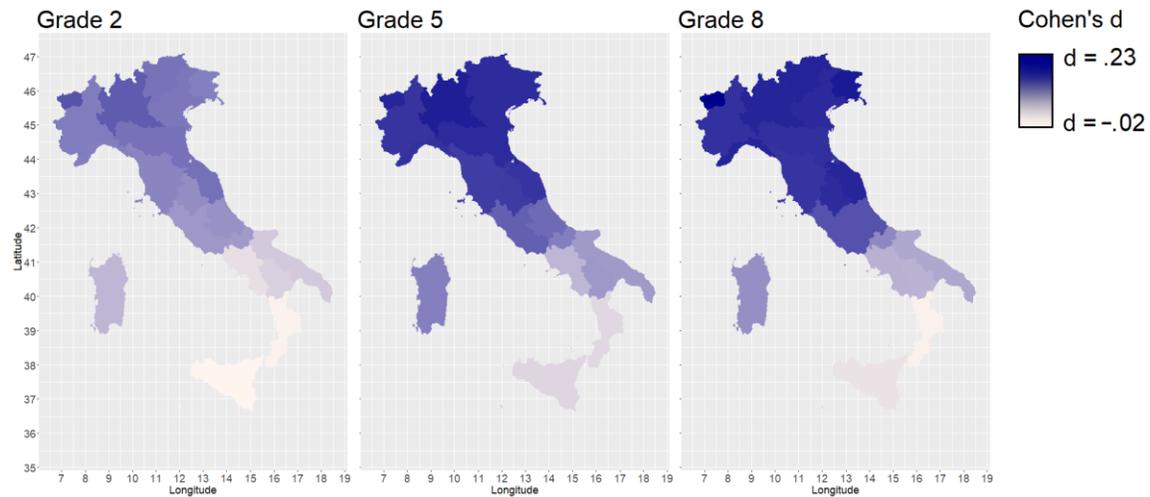

**B.** Males (dis)advantage over females in reading literacy performance

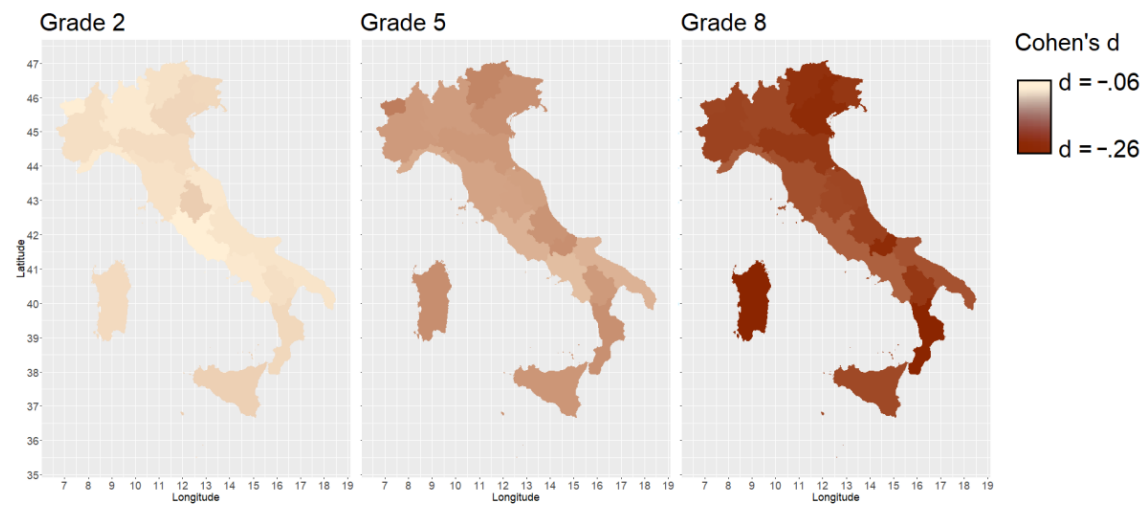

*Figure 1.* Standardized differences (Cohen's d) between males and females in mathematics (A) and reading literacy (B) performance in all Italian regions. Estimates are obtained through meta-analysis over 9 years of INVALSI data on the entire child population of the 2nd, 5th, and 8th grades. Shades of blue indicate males' advantage over females; shades of brown indicate females' advantage over males.



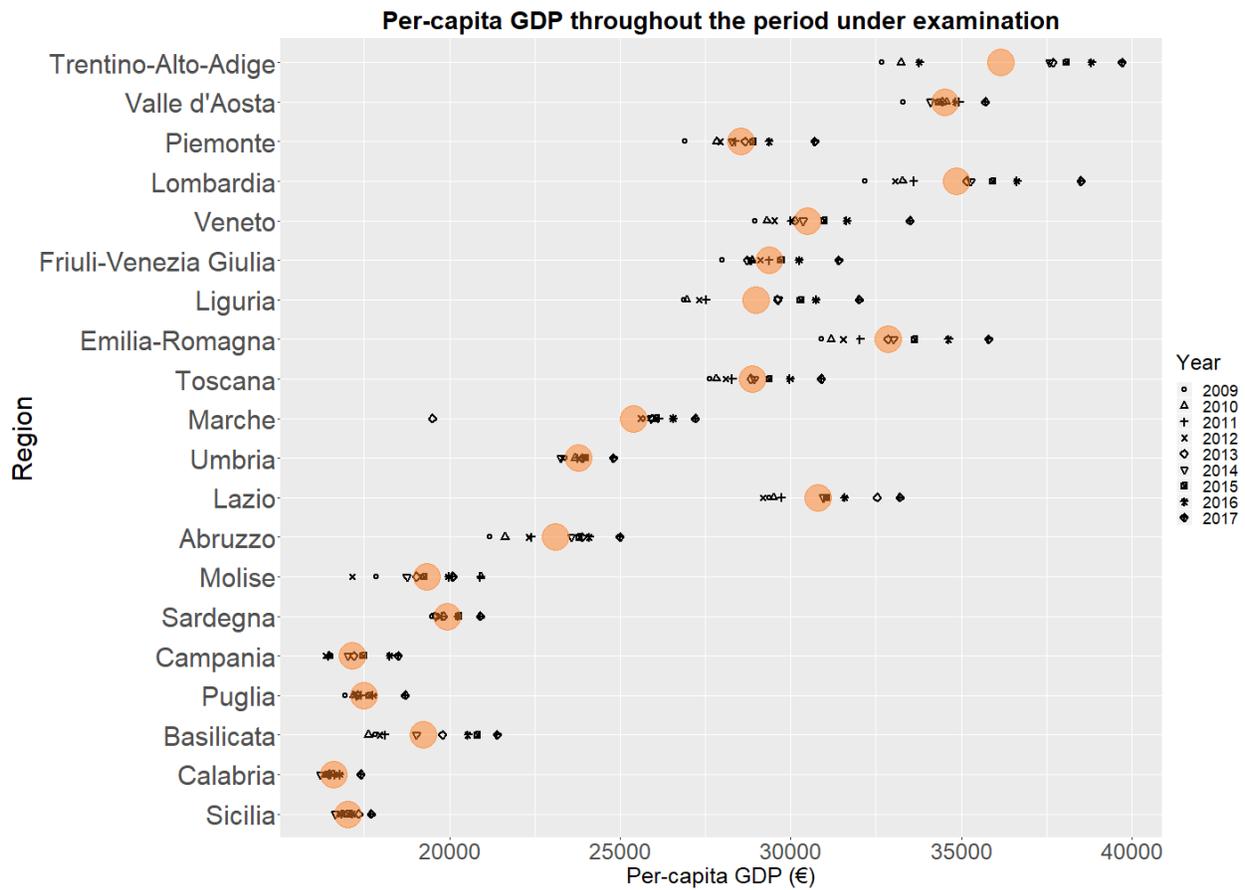

*Figure 2.* Per-capita GPD (€) of all regions throughout the period under examination. The large circular dots represent the grand mean across the years considered. Years considered were 2009-2017 (instead of 2010-2018), because for each year of survey, the level of economic production of the previous year was considered.

**Table 1**. Odds ratios (OR) of Males vs. Females being in the top (or bottom) tails of the distributions of achievement score, as a function of discipline (mathematics, reading) and geographical area (north, center, south), at Grade 8.

| | | Mathematics | | | Reading | | |
|---|---|---|---|---|---|---|---|
| | | North | Center | South | North | Center | South |
| Bottom 1% | OR | 0.83 | 0.79 | 0.89 | 1.64 | 1.44 | 1.66 |
| | 95% CIs | [0.78, 0.89] | [0.73, 0.86] | [0.84, 0.95] | [1.52, 1.77] | [1.30, 1.59] | [1.54, 1.79] |
| | Sim. OR | 0.59 | 0.62 | 0.86 | 1.91 | 1.82 | 1.91 |
| Bottom 5% | OR | 0.77 | 0.76 | 0.88 | 1.59 | 1.49 | 1.58 |
| | 95% CIs | [0.74, 0.81] | [0.72, 0.80] | [0.84, 0.91] | [1.51, 1.67] | [1.39, 1.60] | [1.51, 1.67] |
| | Sim. OR | 0.66 | 0.68 | 0.89 | 1.68 | 1.62 | 1.68 |
| Top 5% | OR | 1.81 | 1.68 | 1.18 | 0.59 | 0.61 | 0.60 |
| | 95% CIs | [1.73, 1.89] | [1.73, 1.89] | [1.13, 1.23] | [0.57, 0.62] | [0.57, 0.64] | [0.58, 0.63] |
| | Sim. OR | 1.52 | 1.46 | 1.13 | 0.59 | 0.62 | 0.59 |
| Top 1% | OR | 2.18 | 1.85 | 1.27 | 0.55 | 0.56 | 0.59 |
| | 95% CIs | [2.04, 2.32] | [1.70, 2.02] | [1.20, 1.36] | [0.51, 0.58] | [0.51, 0.61] | [0.56, 0.63] |
| | Sim. OR | 1.68 | 1.62 | 1.16 | 0.52 | 0.55 | 0.52 |

*Note*. Sim. OR = simulated odd ratio. The simulated odds ratios (OR) are obtained assuming that scores are perfectly normally distributed. The estimates were calculated from 10 million simulated cases that were normally distributed and reproduced exactly the standardized differences between Males and Females observed in the INVALSI data by discipline (mathematics, reading) and geographical area (north, center, south), at Grade 8.

**Supplemental Online Material for the article:**

**Regional effects on the gender gap in mathematics: Results on over 13 million children**

**using the INVALSI dataset**

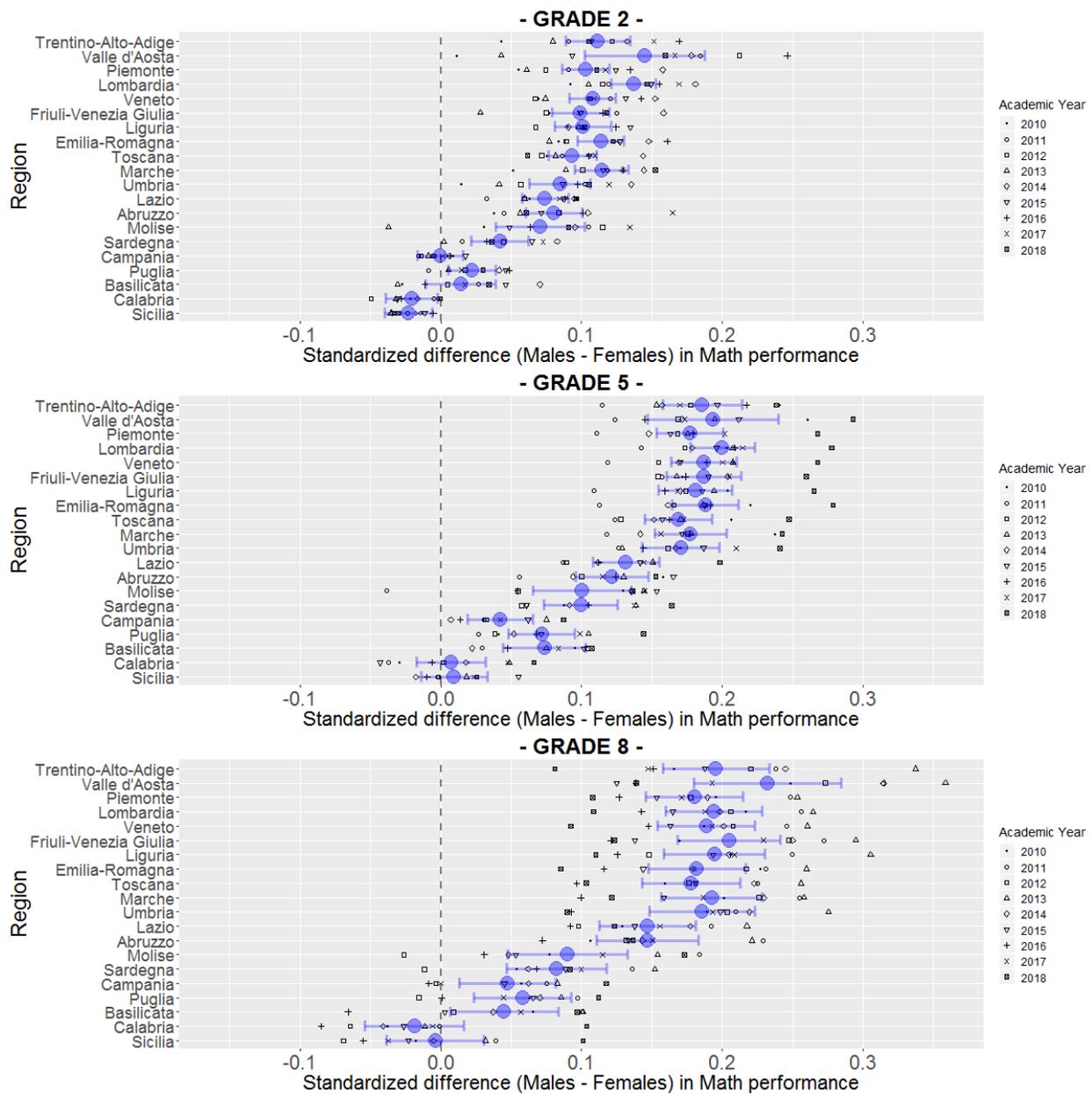

*Figure S1*. Entire population, Mathematics. Meta-analytic estimates of the standardized difference (males vs females) in math performance for each Italian region (divided by grade). These estimates are graphically displayed also in Fig 1 (panel A) of the paper, plotted with different shades of color on the map of Italy. Error bars represent the 95% confidence intervals of the meta-analytic estimates. Smaller points of different shapes represent standardized differences calculated for each year of survey.

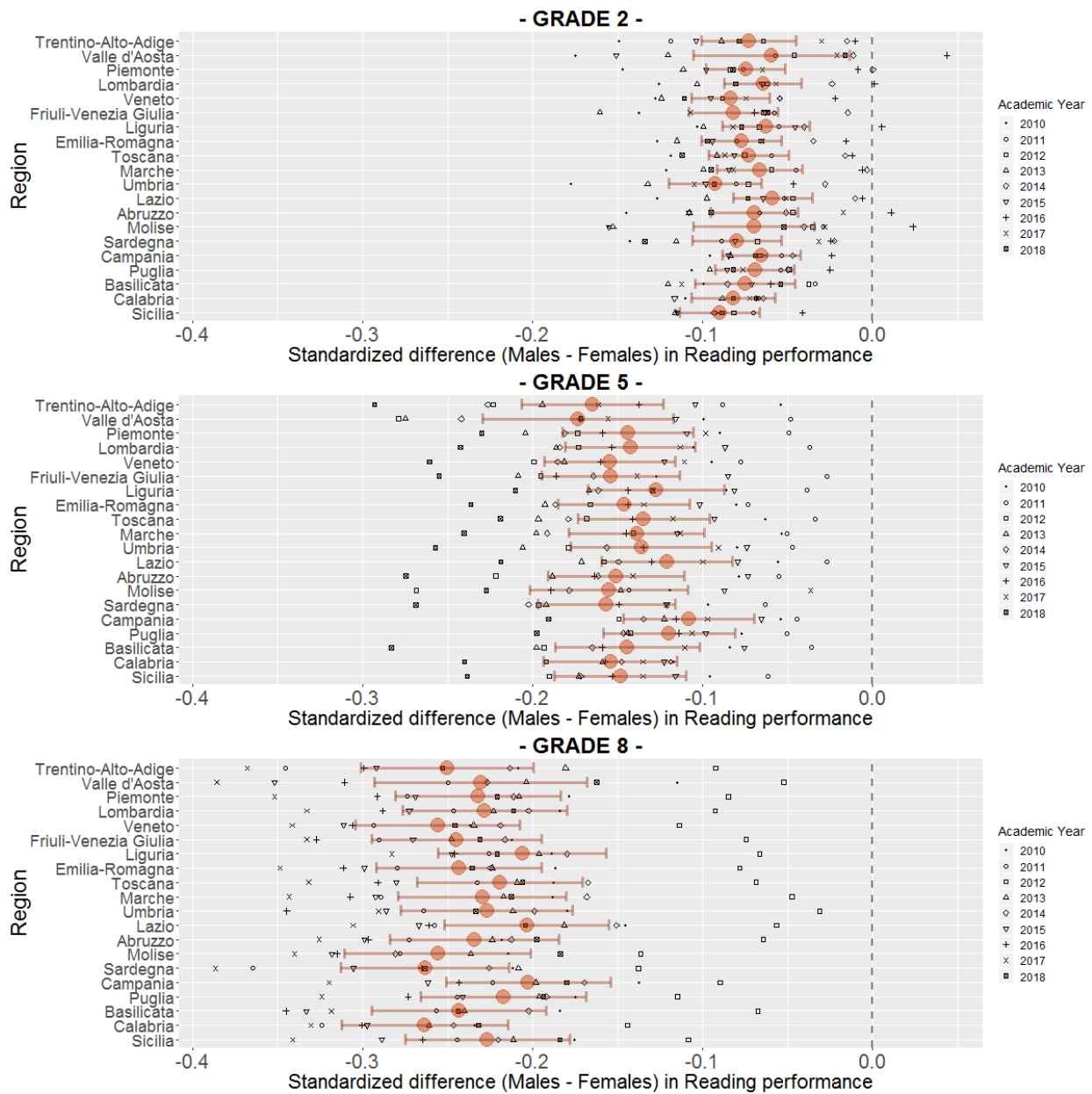

*Figure S2.* Entire population, Reading Literacy. Meta-analytic estimates of the standardized difference (males vs. females) in reading literacy performance for each Italian region (divided by grade). These estimates are graphically displayed also in Fig. 1 (panel B) of the paper, plotted with different shades of color on the map of Italy. Error bars represent the 95% confidence intervals of the meta-analytic estimates. Smaller points of different shapes represent standardized differences calculated for each year of survey.

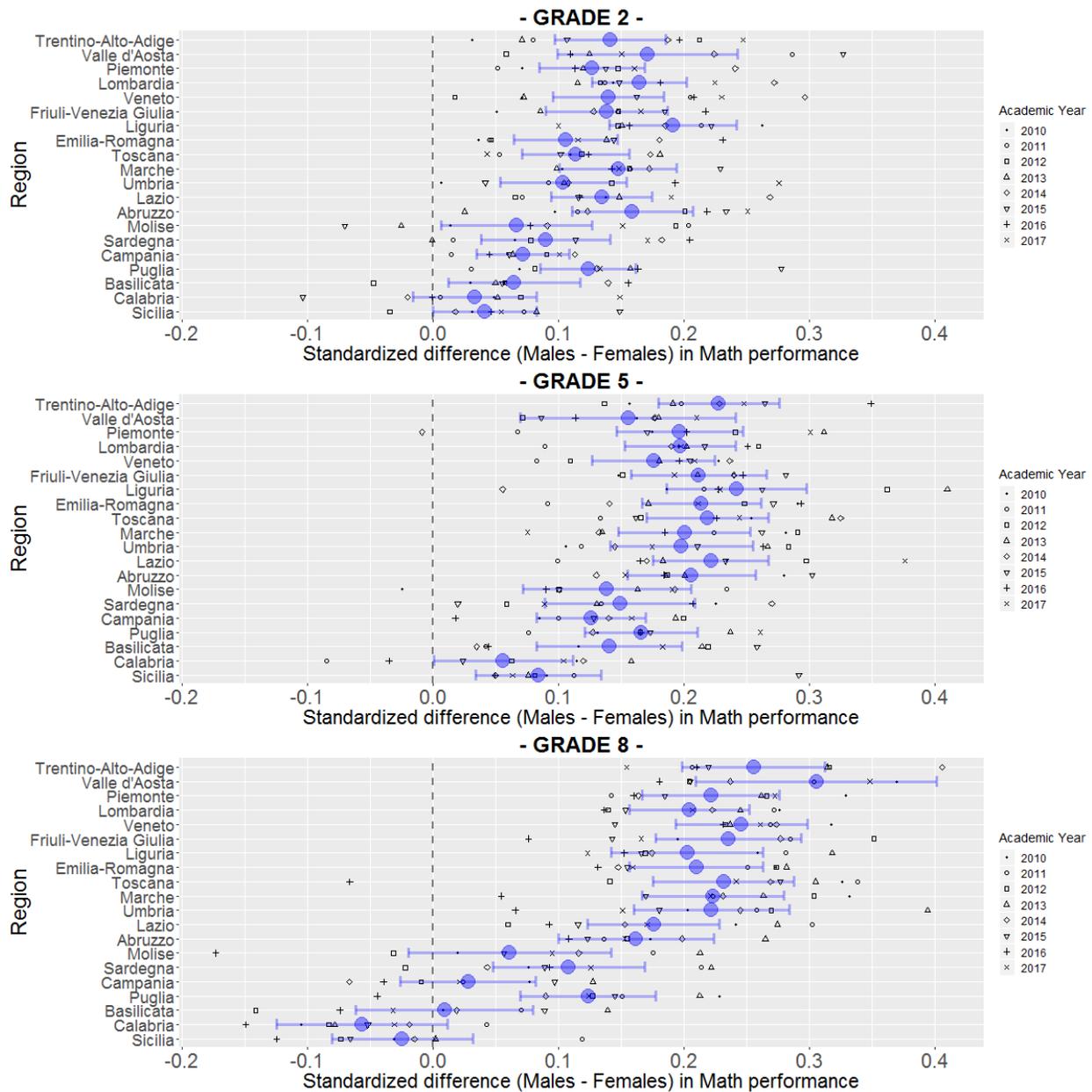

*Figure S3.* Controlled representative sub-sample, Mathematics. Meta-analytic estimates of the standardized difference (males vs. females) in math performance for each Italian region (divided by grade). As is evident, the exact same trends represented in Fig S1 were found on the representative sub-sample, therefore the results gathered from the entire population are robust. In this case, however, the estimates are obviously less precise, as shown by larger 95% confidence intervals (error bars). Smaller points of different shapes represent standardized differences calculated for each year of survey.

**A.** Mathematics

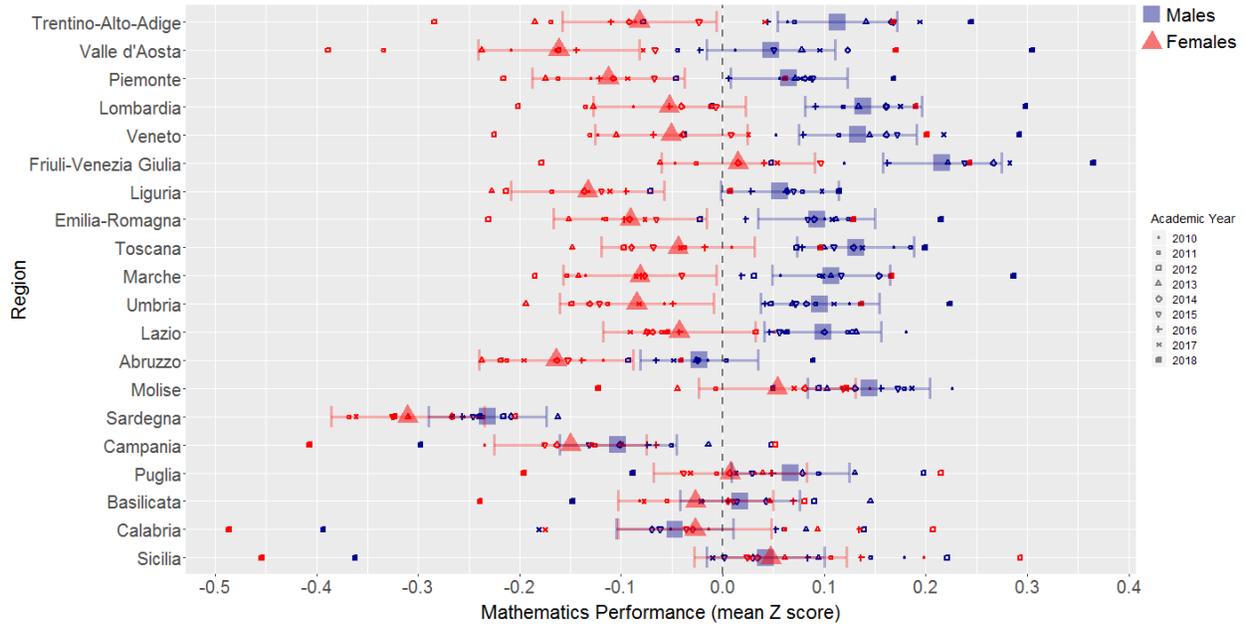

**B.** Reading

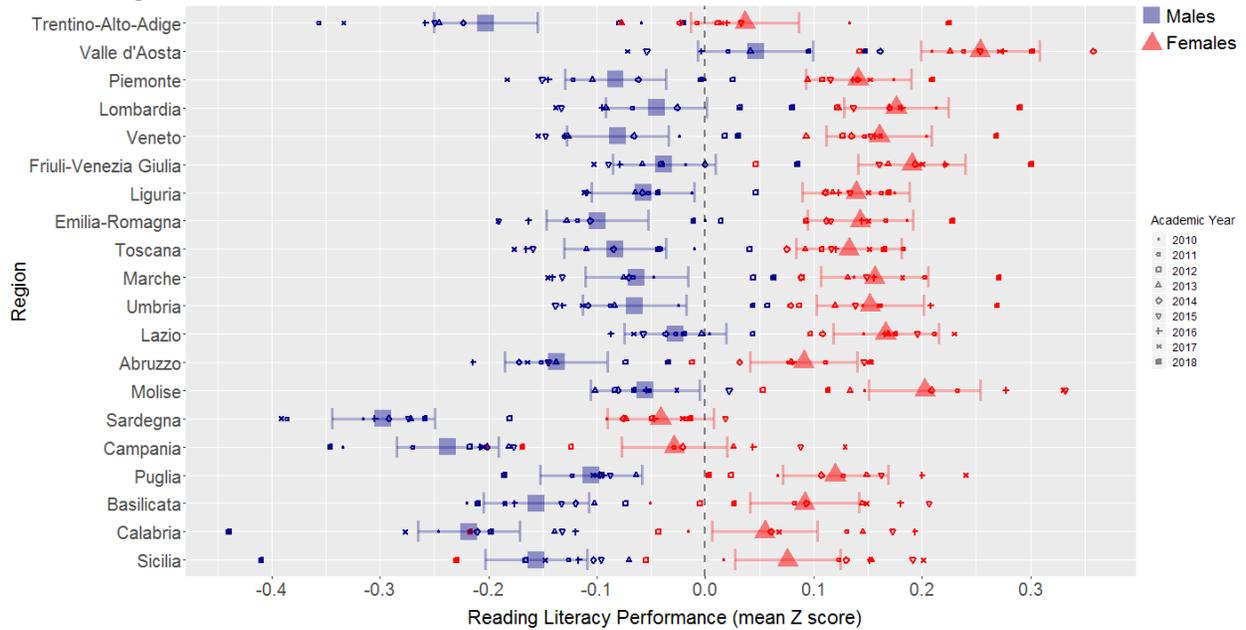

*Figure S4.* Male and female performance in mathematics (A) and reading (B) in the general population in the 8[th] grade. Meta-analytic estimates are based on mean average.

**A.** Employement

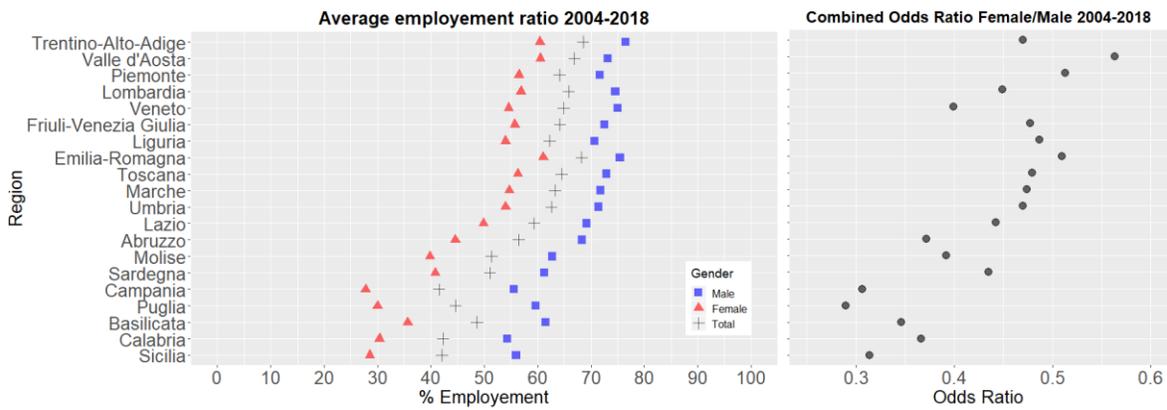

**B.** Unemployement

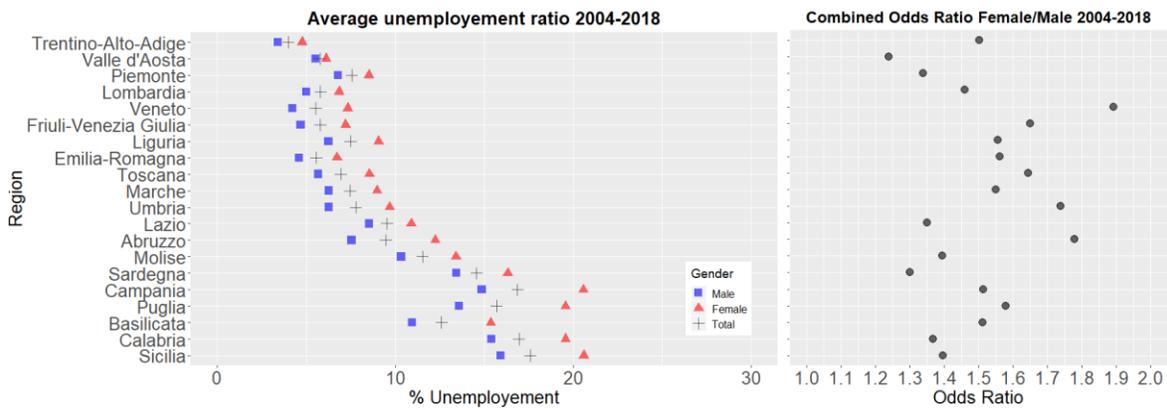

*Figure S5.* Percentage levels of employment (A) and unemployment (B) in all regions, and comparison between males and females. The values refer to the average throughout the period 2004-2018. This is the period that could be retrieved from the ISTAT archive. It corresponds to the period under examination, including data until six years earlier, thus providing a picture of the broader socio-economic background of the period under examination.

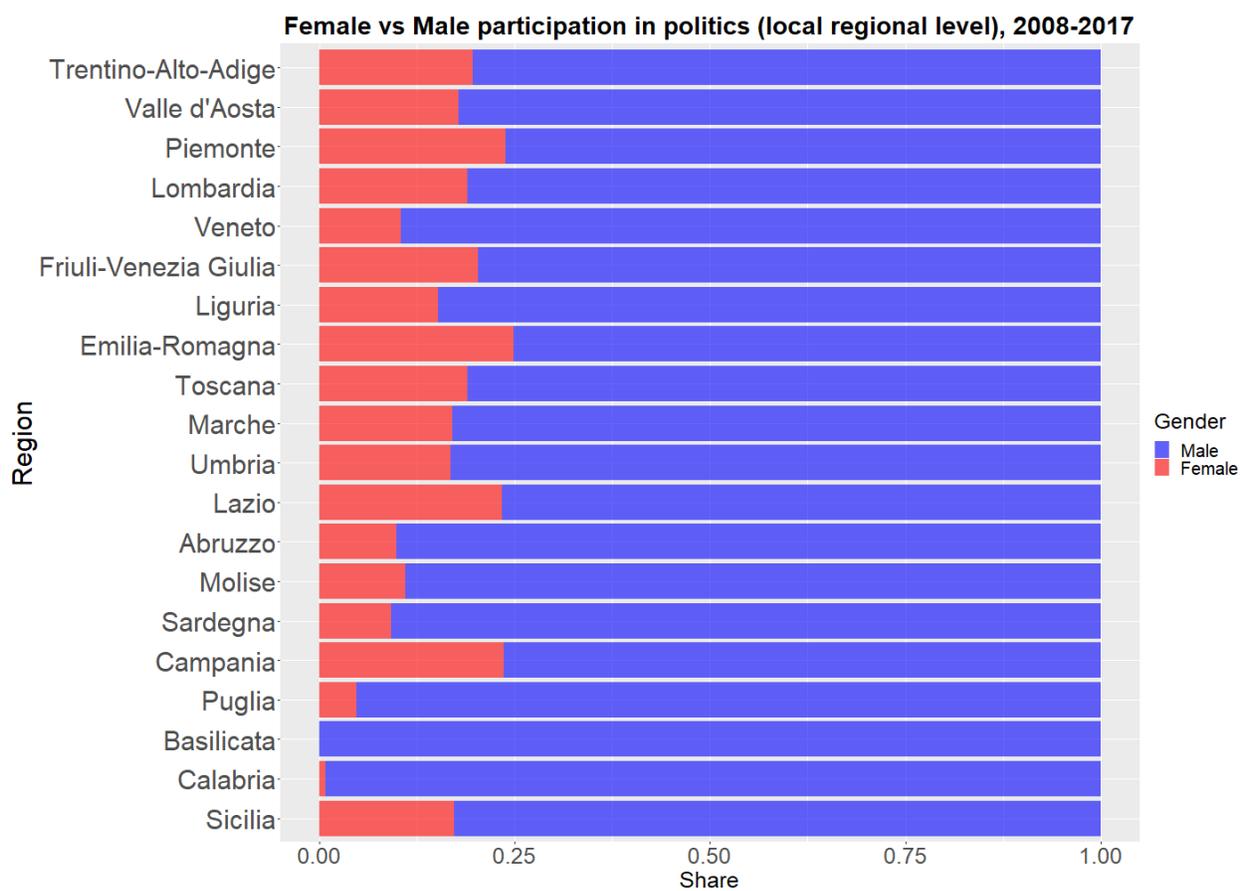

*Figure S6.* Share of female vs male participation in local politics (includes participation in city, province and region Councils). The values refer to the average throughout the period 2008-2017. This is the period that could be retrieved from the ISTAT archive, and it roughly corresponds to the period under examination.

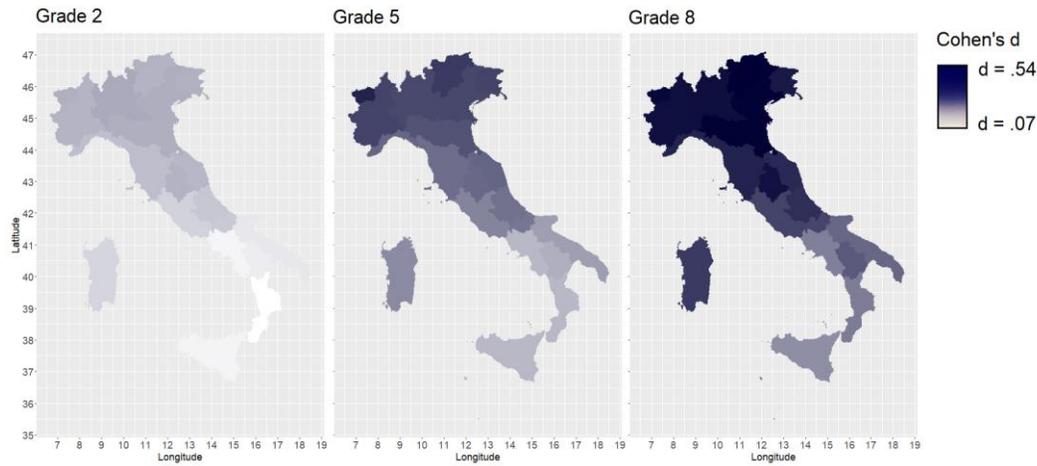

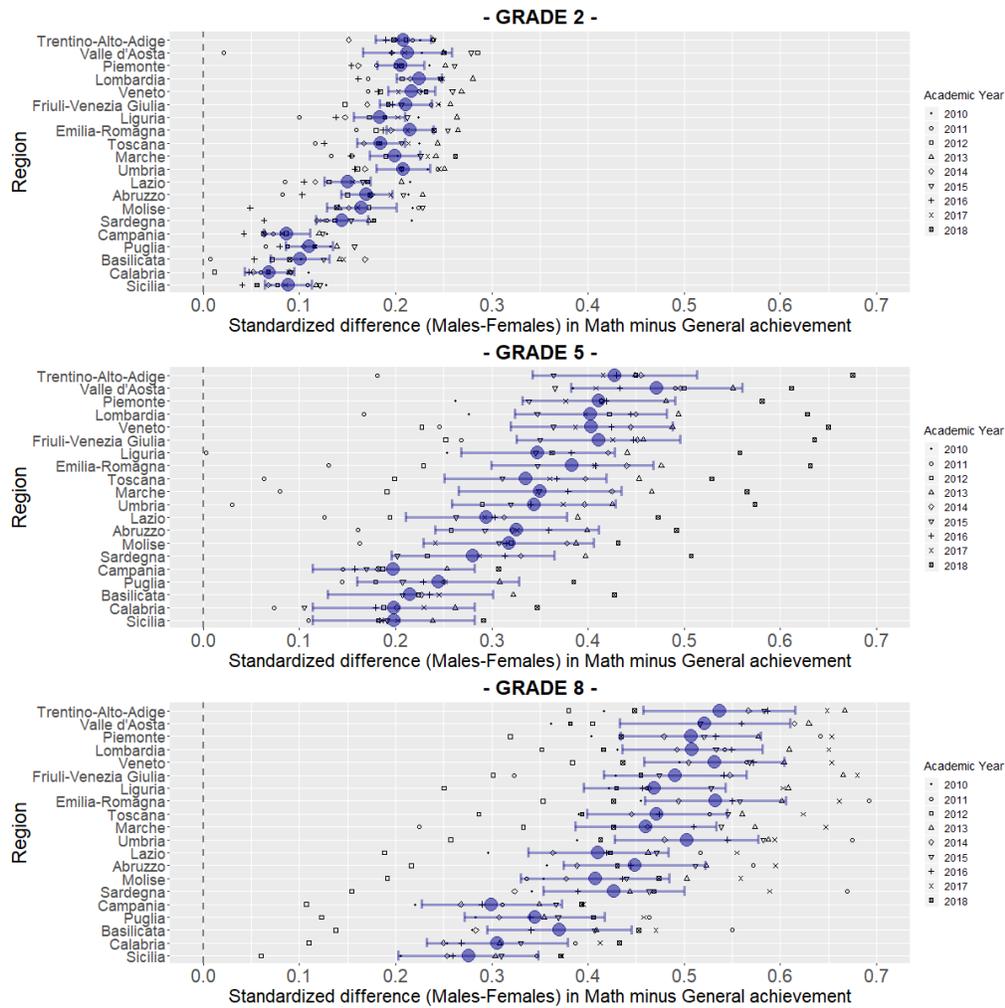

*Figure S7.* Standardized differences between Males and Females in mathematics vis-à-vis overall achievement (i.e., mathematics minus the average between mathematics and reading literacy), with positive scores indicating a male advantage.